\documentclass[11pt,a4paper]{amsart}  

\usepackage[utf8]{inputenc}
\usepackage[T2A]{fontenc}
\usepackage{amsmath, amssymb, euscript, amsthm, amsfonts}
\usepackage{graphicx}
\usepackage{mathrsfs}
\usepackage{latexsym}
\usepackage{tikz-cd}
\usepackage{tikz}
\usepackage{tkz-graph}

\makeatletter
\def\@settitle{\begin{center}%
		\baselineskip14\p@\relax
		\bfseries
		\Large
		\@title
	\end{center}%
}
\makeatother

\usetikzlibrary{automata, positioning, arrows}

\usepackage{epigraph}
\theoremstyle{plain}
\newtheorem{thm}{Theorem}

\newtheorem{prop}{Proposition}

\theoremstyle{definition}
\newtheorem{exer}{Exercise}[section]
\newtheorem{defn}{Definition}
\newtheorem*{ex}{Example}
\newtheorem*{exs}{Examples:}

\theoremstyle{remark}
\newtheorem{rem}{Remark}

\newcommand{\tn}[1]{\textnormal{#1}}
\newcommand{\cat}[1]{\tn{\textbf{#1}}}

\title[]{Quivers, Algebras and Adjoint functors}
\author[]{Kostiantyn Iusenko}
\address{Institute of Mathematics and Statistics\\
	University of São Paulo \\Brazil}
\email{iusenko@ime.usp.br}

\begin{document}

\maketitle 

\setcounter{tocdepth}{1}
\tableofcontents

\section*{About this course}

The aim of this course is a primer  with basic notions of category theory and representation theory of finite dimensional algebras. We will study  the concept of adjoint functors and show the correspondence  ``quiver'' $\leftrightarrows$ ``algebra'' can be interpreted as a pair of adjoint functors between certain categories.

Mostly the lectures do not contain the proofs and the theory is accompanied by examples (sometimes introduced as the exercises). For a deeper acquaintance with the theory of categories and representations of finite dimensional algebras, we recommend (to interested reader) consult with the books and/or on-line resources listed after Lectures 2 and 3.

\section{What are categories and functors?}

\epigraph{``Perhaps the purpose of categorical algebra is to show that which is trivial is trivially trivial''}{--- \textup{Peter Freyd}}

\subsection{Categories.} Most mathematical theories deal with the situations where there are certain maps between objects of a certain nature. The set of objects itself is ``static'', while consideration of morphisms between objects is more ``dynamic''. Usually one imposes the restrictions on the nature of morphisms between objects, for instance, it rarely makes sense to consider all possible maps between groups, usually one limits oneself to studying only group homomorphism.

The concept of a \textit{category} was introduced by Samuel Eilenberg and Saunders Mac Lane as a tool for simultaneous investigation of objects and morphisms between them.
This concept is slightly abstract but very convenient. Before we give the precise definitions, we look at some simple examples.

\newpage

\begin{exs}
\

	\begin{itemize}
		\item Category $ \cat{Sets}$: objects --- set and morphisms --- arbitrary functions between sets;
		\item $ \cat{Groups}$: objects --- groups, morphisms --- homomorphisms of groups;
		\item $ \cat{AbGroups}$: objects --- commutative groups, morphisms --- homomorphisms of groups;
		\item $ \cat{Rings}$: objects --- rings, morphisms --- homomorphisms of rings;
		\item $ \cat{Alg}_k $: objects --- algebras over fixed field $k$, morphisms --- homomorphisms of algebras;
		\item $ \cat{Top}$: objects --- topological spaces, morphisms --- continuous functions;
		\item $ \cat{Mflds}$: objects --- smooth manifolds, morphisms --- differentiable maps between manifolds;
		\item $ \cat{Vect}_k $: objects --- vector spaces over $k$, morphisms --- linear operators.
	\end{itemize}
\end{exs}

It should be noted that in all these examples above we can form a composition of morphisms and this composition is associative (as in all examples morphisms are functions between sets which satisfy certain restrictions and composition of function is associative). 

\begin{defn}
A \textbf{category} $\mathcal C$ consists of the following:
\begin{itemize}
	\item class of objects $\cat{Ob}(\mathcal C)$. The situation ``$X$ is an object in 
	$\mathcal C$'' we write as $X \in \cat{Ob}(\mathcal C)$ or $X \in \mathcal C$;
	\item class of morphisms $\cat{Mor}(\mathcal C)$. Each \textit{morphism} $f$ is a some map from $X \in \mathcal C$ to $ Y \in \mathcal C$. Formally $ \cat{Mor}(\mathcal C)$ is a disjoint union  of classes $ \cat{Mor}(X, Y) $ for all possible $ X, Y \in \mathcal C $. We will denote the morphism $f$ by arrow $X \xrightarrow{\ f\ } Y$. By $ \mathcal C (X, Y) $ we denote the morphisms between $X$ and $Y$.;
	\item composition rule of morphisms:
	\begin{equation*}
	\begin{split}
		\cat{Mor}(X, Y) \times \cat{Mor}(Y, Z) & \rightarrow \cat{Mor}(X, Z), \\
							(f, g) & \mapsto fg,
	\end{split}
	\end{equation*}
	which takes two morphisms $X \xrightarrow{\ f\ } Y$ and $Y \xrightarrow{\ f\ } Z$ to the morphism $X \xrightarrow{\ fg\ } Z$;
	\item for each $ X \in \mathcal C$ there exists an identity morphism $X \rightarrow X $;
\end{itemize}
This structure must satisfy the following axioms:
\begin{itemize}
	\item composition of morphisms is associative;
	\item composition of arbitrary morphism $ f: X \rightarrow Y $ with identity morphism equals $f$.
\end{itemize}
	\end{defn}

\subsection{Functors}

Considering several categories simultaneously, functors bring the way to ``relate'' them.

\begin{defn}
	\textbf{Covariant} (resp. \textbf{contravariant}) \textbf{functor} $F$ from a category $\mathcal C $ to a category $\mathcal D$ is a rule that associates to an arbitrary object $ X \in \mathcal C $ an object $ F (X) \in \mathcal D $, and to an arbitrary morphism $ f: X \rightarrow Y $ a morphism $ F (f): F (X) \rightarrow F (Y) $ (resp. $ F (f): F (Y) \rightarrow F (X) $) such that the following axioms hold:
\begin{itemize}
	\item $ F (\tn{id}_X) =\tn{id}_{F (X) }$ for any $ X \in \mathcal C $;
	\item $ F $ preserves the composition between morphisms, that is, for arbitrary $ f: X \rightarrow Y $, $ g: Y \rightarrow Z $ are $ F (f \circ g) = F ( f) \circ F (g) $ (if $ F $ covariant) and $ F (f \circ g) = F (g) \circ F (f) $ (if $ F $ contravariant) .
\end{itemize}	
\end{defn}

\begin{exs}
\

\begin{itemize}
	\item
		\textbf{inclusion of subcategory.} Let $ \mathcal C $  be a subcategory  of $ \mathcal D $. That means that objects and morphisms in $ \mathcal C $ are also objects and morphisms in $ \mathcal D $. The inclusion functor $ \mathcal C \hookrightarrow \mathcal D $, acts identically on objects and morphisms.
		For example, we have a functor
		$$
			\cat{AbGroups}\hookrightarrow \cat{Groups}.
		$$
	\item \textbf{forgetful functor:} define the functor 
	$ F: \cat{Groups}\rightarrow \cat{Set}$, which maps an arbitrary group to its underlying set (forgetting group structure), and arbitrary homomorphism between groups to a function between the respective sets. Functor ``forgets'' group structure (on the objects and morphisms). Similarly, we have the forgetful functor $ \cat{Vect}_k \rightarrow \cat{Sets}$, which maps any vector space to the set of all its vectors, and any linear transformation between the spaces to the corresponding map between the sets of vectors. In fact, we forget that we can add vectors and multiply them by scalars, and that linear maps are linear.
Similarly, one can define the forgetful functors: $ \cat{Alg}_k \to \cat{Set}$, $ \cat{Top}\to \cat{Set}$, $ \cat{Ring}\to \cat{Set}$, $ \cat{R-Mod}\to \cat{Set}$, ..., $ \cat{Alg}_k \to \cat{AbGroup}$ ( functor which forgets multiplication).
	
	\item \textbf{Free functors:} For any set $X$ we define $F(X)$ as a free group generated by the set $X$. Arbitrary function $f:X \to Y$, which maps $x \in X $ to $f (x) \in Y$, defines the group homomorphism $F(f): F (X)  \to F (Y) $. It is easy to see that this map satisfies $F(fg) = F (f)F(g)$, hence it defines the functor $F: \cat{Sets}\to \cat{Groups}$.
Similarly, one defines other ``free'' functors: $ \cat{Sets}\to \cat{Vect}_k $, 
	$\cat{Sets}\to \cat{Ring}$, ..., $ \cat{Sets}\to \cat{Top}$ (set $X$ endowed with the discrete topology), and many others.
	\item \textbf{functor Boolean}. We define the functor $ \mathcal P: \cat{Sets}\to \cat{Sets}$, setting $ \mathcal P [X] $ as the set of all subsets of $X$. Now if $ f: X \to Y $ is a function between sets, and $U \subset X$, define $ \mathcal P [f](U)$ as the image of $U$ under $f$.
	\item Consider the following example of contravariant functor. Duality functor
	$$ \cat{Vect}_{k}\rightarrow \cat{Vect}_{k}, $$ which maps an arbitrary vector space $ V $ to vector space $ V ^ * $ of all linear functional on $ V $. Linear operator $ L: U \rightarrow V $ is mapped  to its conjugate $ L ^ *: V ^{* }\rightarrow U ^{* }$ (which maps arbitrary linear functional $ \varphi \in V ^ * $ to the linear functional $ u \mapsto \varphi (L (u)) $ on $ U $).
	\item Another example of contravariant functor is the functor $ \cat{Top}\rightarrow \cat{Rings}$, which assigns to each topological space $ X $ the ring of continuous functions $ C ^ 0 (X, \mathbb R) $, and to an arbitrary continuous map $ f: X \rightarrow Y $ it assigns so-called (pull-back) map $ f ^ *: C ^ 0 (Y, \mathbb R ) \rightarrow C ^ 0 (X, \mathbb R) $ (observe that the composition of function on $Y$ with $ f $ is a function on $X$).

\end{itemize}
\end{exs}

\subsection{Equivalence of categories}

Given two categories $ \mathcal C $ and $ \mathcal D $ when are they equivalent? It is natural to say that $ \mathcal C $ and $ \mathcal D $ are isomorphic if there are functors $ F: \mathcal C \rightarrow \mathcal D $ and $ G: \mathcal D \rightarrow \mathcal C $ such that they are inverses to each other. In fact, this definition is quite restrictive, as the following example explains:

\begin{ex}\label{linalg-cat}
Let $\mathcal D$ be the category of finite dimensional vector spaces over $k$, and $ \mathcal C $ be a subcategory containing the vector space $ k ^ n $ (of column vectors) for each dimension $ n $. Note that the $ \cat{Mor}(k ^ n, k ^ m) $ can be naturally identified with matrices $M_{m,n}(k)$. The categories $ \mathcal C $ and $ \mathcal D $ are not isomorphic as $ \mathcal D $ contains all possible vector spaces. But an arbitrary $n$-dimensional vector space $ V $ is isomorphic to $ k ^ n $ (selecting the basis in $V$), so the category of $ \mathcal C $ in some sense full enough and we can consider $ \mathcal C $ and $ \mathcal D $ as equivalent categories.
\end{ex}

To formalize the last example, consider the following definition

\begin{defn}\label{cat-equiv}
	Covariant functor $ F: \mathcal C \rightarrow D $ is called the \textbf{equivalence of categories} if
	\begin{itemize}
		\item $ F $ is \textbf{essentially surjective}, that is an arbitrary object in $ \mathcal D $ is isomorphic (but not necessarily equal!) to an object of the form $F(X)$ for some $X \in \mathcal C$.  
		\item $F$ is \textbf{full} and \textbf{faithful}, that is there is a bijection
		$$
			\cat{Mor}_{\mathcal C}(X, Y) \simeq \cat{Mor}_{\mathcal D}(F (X), F (Y)),
		$$
		for arbitrary $ X, Y \in \mathcal C $.
	\end{itemize}
\end{defn}

Consider the category as in Example \ref{linalg-cat}. We make sure that the inclusion functor $ F: \mathcal C \rightarrow \mathcal D$ is an equivalence of categories.
Indeed, if $V$ is an arbitrary $n$-dimensional space, $ V $ is isomorphic to $k ^ n$. Fixing a basis $ e_1, \dots, e_n $, we get the isomorphism $ V \rightarrow k ^ n $, which maps an arbitrary vector $ v \in V $ to the column of coordinates of $v$ in the basis $ \{e_i \}$, and therefore $F$ is essentially surjective. It is easy to see that $ F $ is full and faithful, that is $ F $ is an equivalence of categories.

\subsection{Exercises}

\begin{exer}
	Let $ \mathcal C $ be an arbitrary category. Show the following:
	\begin{itemize}
		\item [a)] identical morphism $ id_X: X \rightarrow X $ is unique for each object $ X \in \mathcal C $;
		\item [b)] an arbitrary isomorphism $ \mathcal C $ has unique inverse;
		\item [c)] Let $ f: X \rightarrow Y $ and $ f: Y \rightarrow Z $ be two morphisms. Show that if two of the morphisms $ f, g $ and $ f \circ g $ are isomorphisms, then the third is also isomorphism. (This property is called \textit{two of three}).
		\ 
	\end{itemize}
\end{exer}

\begin{exer}\label{posetCat}
Let $I$ be an arbitrary partially ordered set in which a partial order is given by $\preceq$. With $I$ we associate the category $ \mathcal C_I $, in which objects are elements of the set $I$, and for two arbitrary $i, j \in I$, $\cat{Mor}(i, j)$ --- empty set if $ i \npreceq j $ and has one element if $i \preceq j$.
Using reflexivity and transitivity of relation $\preceq$, one determines the composition of morphisms in $ \mathcal C_I $ and shows that the construction above defines the category $ \mathcal C_I$.
\end{exer}

\begin{exer}\label{posetFuct}
Let $I$ and $J$ be two posets. Show that an arbitrary functor between the categories $ \mathcal C_I $ and $ \mathcal C_J $ is given by poset homomorphism (i.e. by the function that preserves order) $f: I \to J$.	
\end{exer}

\begin{exer}\label{posetTop}
Let $X$ be an arbitrary topological space and $I(X)$ is the set of all closed subsets of $X$. Show that $I(X)$ is a poset (with the order given by the inclusion of subsets). Specify category $ \cat{Top}(X) $ as $ \mathcal C_{I (X) }$ from Exercise \ref{posetCat}.
\end{exer}

\begin{exer}\label{complexification}
	Let $ V $ be vector space over $\mathbb R$. Show that the \textit{complexification} 
	$ V \mapsto V \otimes_{\mathbb R}\mathbb C $ defines the functor $ \cat{Vect}_{\mathbb R}\to \cat{Vect}_{\mathbb C}$. 
\end{exer}

\begin{exer}
	Let $ \mathcal C, \mathcal D $ be arbitrary categories, and $ F: \mathcal C \rightarrow \mathcal D $ be a functor. Show the following:
	\begin{itemize}
		\item [a)] $ F $ maps isomorphisms to isomorphisms;
		\item [b)] if $ F $ is full and faithful, and $ F (f): F (X) \rightarrow F (Y) $ is an isomorphism in $ \mathcal D $, then $ f: X \rightarrow Y $ is an isomorphism in $\mathcal C$.
	\end{itemize}
\end{exer}

\begin{exer}\label{functorAb}
	Let $ G $ be an arbitrary group and $ [G, G] $ its commutator. Show that $ G / [G, G] $ is an abelian group, and the map $ G \mapsto G / [G, G] $ defines the functor $ Ab: \cat{Group}\to \cat{AbGroup}$.
\end{exer}

\begin{exer}
Define $ \cat{Top}_{* }$ as a category whose objects are the pairs $ (X, x_0) $, $ X $ is topological space, $ x_0 \in X $ is a fixed point and morphisms $ f: (X, x_0) \to (Y, y_0) $ are continuous maps $ f: X \to Y $ such that $ f (x_0) = y_0 $. Make sure that $ \cat{Top}_{* }$ is a category, and the map $ \pi_1: \cat{Top}_{* }\to \cat{Group}$, which to the pair $ (X, x_0) $ connects its fundamental group $ \pi_1 (X, x_0) $, defines covariant functor.
\end{exer}

\begin{exer}
Define $ \cat{CHaus}$ a category whose objects are compact Hausdorff spaces and morphisms are continuous maps. And define the map $ \beta: \cat{Top}\to \cat{CHaus}$, which to an arbitrary topological space $X$ associates its Stone-Cech compactification $ \beta X $ (i.e., ``maximum'' compact Hausdorff space ``generated'' by $X$). Show that $\beta$ defines a functor.
\end{exer}

\begin{exer}\label{homten}
Let $ R, S $ be rings (not necessarily commutative). Consider the category of right modules over these rings $ \mathcal C =\textrm{Mod}_R $ and $ \mathcal D =\textrm{Mod}_S $ (in which morphisms are homomorphisms of modules). Show that, fixing $(R, S)$-bimodule $X$, one can define two functors $ F: \mathcal C \to \mathcal D $ and $ G: \mathcal D  \to \mathcal C $ as follows:
	\begin{align*}
		F (Y) & = Y \otimes_R X, \qquad Y \in \mathcal D, \\
		G (Z) & =\textrm{Hom}_S (X, Z), \qquad Z \in \mathcal C.
	\end{align*}
\end{exer}

\newpage

\section{Adjunction between categories}

\subsection{Natural transformations}

\epigraph{`` I did not invent category theory to talk about functors. I invented it to talk
about natural transformations. '' }{--- \textup{Saunders Mac Lane}}

What is a natural transformation? This is a map from one functor to another! Consider the simple example that explains this quite vague definition.
Recall that for any vector space $V$ there exists (``natural'') linear map
\begin{align*}
	\phi_V: & V \to V ^{** }\\
			& V \mapsto (f \mapsto f (v)).
\end{align*}
If $ \dim V <\infty $, it is easy to see that $\phi_V$ an isomorphism. What naturalness of $\phi_V$ means? In fact, if $\dim V<\infty $, then the spaces $V$ and $V^*$ are isomorphic, but there is no canonical isomorphism (the isomorphisms requires the choice of basis in  $V$). On the other hand, the isomorphism $\phi_V$ does not require any additional choice. To formalize this construction, consider the properties of $\phi_V$. An arbitrary linear map $L: U \to V$ generates a linear map $L ^{** }:U ^{** }\to V ^{** }$, which, together with the maps $ \phi_U, \phi_V $ generates the following diagram:
\begin{center}
\begin{tikzcd}[row sep = huge, column sep = huge]
U \arrow [r, "\phi_U"] \arrow [d, "L" '] & U ^{** }\arrow [d, "L ^{** }"] \\
V \arrow [r, "\phi_V"] & V ^{** }
\end{tikzcd}
\end{center} 
There is no a priori reason for supposing that the diagram above is commutative (if $\phi_U$ is an arbitrary linear map, it is clear that the diagram is not commutative). However, the diagram is indeed commutative! Let us show that. Let $u \in U$, we check that
\begin{equation}\label{cd-natural}
	\phi_V (L (u)) = L ^{** }(\phi_U (u)).
\end{equation}
Indeed, the functional in left side of \eqref{cd-natural} maps arbitrary linear functional $f \in V^{*}$ to the value $ f (L (u)) $. While the functional from the right side takes an arbitrary linear functional $f \in V^{*}$ to
$$
	\phi_U (u) (L ^{*}(f)) = L^{*}(f)(u) = f(L(u)).
$$
We encourage the readers to analyse these equalities on their own. 

\begin{defn}
Let $ F, G: \mathcal C \to \mathcal D $ be two covariant functors. \textbf{Natural transformation} $\alpha: F \to G $ between them is a rule that to every object $ X \in \mathcal C $ associates a morphism $\alpha_X: F (X) \to G (X) $ such that for every morphism $ f: X_1 \to X_2 $ in $\mathcal C$ the following diagram commutes
\begin{center}
\begin{tikzcd}[row sep = huge, column sep = huge]
F (X_1) \arrow [d, "F (f)" '] \arrow [r, "\alpha_{X_1}"] & G (X_1) \arrow [d, "G (f ) "] \\
F (X_2) \arrow [r, "\alpha_{X_2}"] & G (X_2)
\end{tikzcd}
\end{center} 
If $ \alpha_X $ is an isomorphism for every object $X$, then $\alpha$ is called natural \textbf{isomorfism}. 
\end{defn}

Let's analyse the previous example in the context of this definition. Consider the category of vector spaces over a field $k$ denoted by $\cat{Vect}_k$ and two functors: $ \textrm{Id}: \cat{Vect}_k \to \cat{Vect}_k $ (identical) and $D: \cat{Vect}_k \to \cat{Vect}_k $ (double duality), which maps an arbitrary vector space $V$ to its double dual $V^{**}$, and arbitrary morphism $L: U \to V$ to $ L ^{** }: U ^{** }\to V ^{** }$.
\begin{exer}
Show that $\phi_U$ defines natural transformation between functors $\textrm{Id}$ and $D$, and considering subcategory $\cat{FVect}_k \subset \cat{Vect}_k$ of finite dimensional vector spaces $\phi_U $ defines a natural isomorphism between the respective restricted functors. 
\end{exer}

\subsection{Adjoint functors}

Recall that there are two functors between the categories $ \cat{Sets}$ and $ \cat{Vect}_k $: free $ F: \cat{Sets}\to \cat{Vect}_k$ and forgetful: $\cat{Vect}_k \to \cat{Sets}$. $ G (V) $ is the set of all vectors for a given vector space $V$, and $F(X)$ is a vector space over the field $k$ with basis $X$ (i.e. $F(X)$ is formed by all formal linear combinations $\sum_i \lambda_i x_i $ with $\lambda_i \in k$ and $x_i \in X$, endowed with obvious vector space structure). An arbitrary function $g: X \to G (V)$ can be uniquely extended to a linear operator $f: F (X) \to V$ (the operator $f$ defined as follows $f(\sum_i \lambda_i x_i) =\sum_i \lambda_i g (x_i)$). This generates the map $\eta: g \mapsto f$ which has an ``inverse'' $\mu: f \mapsto f |_X$ (which assigns to an arbitrary linear operator $f: F(X) \to V$ the map $f|_X: X \to G (V) $, restricting $f$ to the basis of $X$). Thus, $\eta =\eta_{X, V}$ defines the bijection
\begin{equation}\label{eta}
	\eta: \cat{Sets}(X, G (V)) \cong \cat{Vect}_k (F (X), V)
\end{equation}
Moreover, this bijection defined in ``canonical'' way for all sets $X$ and vector spaces $V$, i.e. component-wise $\eta_{X, V}$ defines the natural transformation of functors if we regard left and right sides of $\eqref{eta}$ as the functors  in variables $X$ and $V$. A more detailed interpretation we give in definition of \textit{adjunction}.

\begin{defn}
	Suppose we are given two functors $ F: \mathcal C \to \mathcal D $ and $ G: \mathcal D \to \mathcal C $. An \textbf{adjunction} between $ F $ and $G$ is a law which for any pair $(A \in \mathcal C, B \in \mathcal D)$ associated  the bijection $\eta_{A, B}$ between  $\mathcal C (A, G (B))$ and $\mathcal D(F(A),B)$, which is natural in $A$ and in $B$. In this case the functor $F$ is called \textbf{left adjoint} to $G$, and $G$ functor called \textbf{right adjoint} to $F$.
\end{defn}

It is easy to see that naturalness of bijection $\eta$ means that for each $f: A \to A'$ and $g: B \to B'$ the following diagrams commute:

\begin{center}
\begin{tikzcd}[row sep = huge, column sep = huge]
\mathcal C (A, G (B)) \arrow [d, "f ^ *" '] \arrow [r, "\eta_{A, B}"] & \mathcal D (F (A), B) \arrow [d, "(F (f)) ^ *"] \\
\mathcal C (A ', G (B)) \arrow [r, "\eta_{A', B}" '] & \mathcal D (F (A'), B)
\end{tikzcd},
\end{center} 
and
\begin{center}
\begin{tikzcd}[row sep = huge, column sep = huge]
\mathcal C (A, G (B)) \arrow [d, "(G (g)) _ *" '] \arrow [r, "\eta_{A, B}"] & \mathcal D (F (A), B) \arrow [d, "g_ *"] \\
\mathcal C (A, G (B ')) \arrow [r, "\eta_{A, B'}"'] & \mathcal D (F (A), B ')
\end{tikzcd}.
\end{center} 
Here we used short notation $g_*=\mathcal D(F(A), g)$ (composition with morphism $g$), and $f^* =\mathcal C(f, G (B))$ (pre-composition with $f$).

\begin{rem}
The name ``adjoint functor'' arises as a kind of generalization of the concept of 
``adjoint operator''. Indeed, consider the set of morphisms as bi-functor
	$$
		\cat{Mor}: \mathcal C ^{op}\times \mathcal C \to \cat{Sets}
	$$
compare with the definition of adjoint functors with the definition of adjoint operator in a complex vector space $V$ with scalar product
   	$$
   		\langle \cdot, \cdot \rangle: V ^ c \times V \to \mathbb C,
   	$$
in which $V^c$ denotes the vector space in which the action of the field $\mathbb C$
is precomposed with complex conjugation.
\end{rem}

\begin{exs}
\

\begin{itemize}
	\item [a)] The pairs of functors ``free--forgetful'' are the typical examples of adjunction between functors. Free functor is left adjoint to forgetful and forgetful is  right adjoint to free. One example we saw at the beginning of this section. A similar example: free functor $ F: \cat{Sets}\to \cat{Groups}$ is left adjoint to forgetful functor $ G: \cat{Group}\to \cat{Sets}$. Indeed, any group homomorphism $F(X) \to U$ uniquely and naturally (check!) gives rise to the function $X \to G (U) $. Similarly, we can build examples of adjoint functors for free rings, free $R$-modules and so on.
	\item [b)] forgetful functor $G: \cat{Top}\to \cat{Sets}$ has a left and right adjoint functors. Left adjoint functor $L$ endows  the set $X$  with  discrete topology (because all maps $L(X) \to Y $ are continuous, for an arbitrary $Y \in \cat{Top}$). Right adjoint functor $R$ endows  the set $X$ with trivial topology.
	\item [c)] Let $I$ and $J$ be two posets. Any functor $ F: \mathcal C_I \to \mathcal C_J $ is a an order preserving map $F:I\to J$ (see. Exercise \ref{posetFuct}). Therefore a pair of adjoint functors $F: \mathcal C_I \to \mathcal C_J $, $ G: \mathcal C_J \to \mathcal C_I$ is a pair of order-preserving maps satisfying
	$$
		\mathcal C_I (F (a), b) \cong \mathcal C_J (a, G (b))
	$$
for all $ a \in I $ and $ b \in J $. On the other hand, it means that
$$
	F (a) \preceq_J b \Leftrightarrow a \preceq_I G (b).
$$
The latter correspondence is called Galois correspondence between posets 
(see. Exercises at the end this lecture for specific examples of such correspondence).
\item [d)] functor $ Ab: \cat{Group}\to \cat{AbGroup}$ (see. Exercise \ref{functorAb}) is left adjoint to embedding functor $G:\cat{AbGroup} \to \cat{Group}$.

\item [f)] functor Stone-Cech compactification $ \beta: \cat{Top}\to \cat{CHaus}$ is left adjoint to inclusion functor $ \cat{CHaus}\hookrightarrow \cat{Top}$.

\item [g)] 
\textbf{Tensor-Hom adjunction.} 
Let $U, V, W$ be three vector spaces. Standard fact from linear algebra says that there is an isomorphism
$$
	\textrm{Hom}(U \otimes V, W) \cong \textrm{Hom}(U, \textrm{Hom}(V, W)),
$$
where $ \textrm{Hom}(\cdot, \cdot) $ denotes the vector space of all linear operators. 
Similar correspondence (an isomorphism of abelian groups) holds for right modules over rings (see. Exercise \ref{homten}):
$$
	\textrm{Hom}_S (Y \otimes_R X, Z) \cong \textrm{Hom}_R (Y, \textrm{Hom}_S (X, Z)).
$$
Hence the functor $-\otimes_R X$ (see. Exercise \ref{homten}) is left adjoint to functor 
$\textrm{Hom}_S(X, -)$. 
\end{itemize}	

\end{exs}

The following theorem holds.
\begin{thm} 
To specify the adjunction between two functors $ F: \mathcal C \to \mathcal D $ and $G: \mathcal D \to \mathcal C $ is equivalent to specify two natural transformations $\eta: \textrm{Id}_{\mathcal C}\to GF $ and $ \epsilon: FG \to \textrm{Id}_{\mathcal D}$ such that the following diagrams commute:
\begin{center}
\begin{tikzcd}[row sep = huge, column sep = huge]
F \arrow [transform canvas ={xshift = 0.3ex}, -]{dr}\arrow [transform canvas ={xshift = -0.4ex}, -]{dr}\arrow [r, "F \eta "] & FGF \arrow [d," \epsilon F "] \\
 & F
\end{tikzcd}, and \\
\begin{tikzcd}[row sep = huge, column sep = huge]
G \arrow [transform canvas ={xshift = 0.3ex}, -]{dr}\arrow [transform canvas ={xshift = -0.4ex}, -]{dr}\arrow [r, "\eta G "] & GFG \arrow [d," G \epsilon "] \\
 & G
\end{tikzcd}.
\end{center}
\end{thm}
The natural transformation $\eta$ is called \textbf{unit} of adjunction and $\epsilon$  is \textbf{counit} of adjunction. 

\begin{ex}
Consider adjoint pair of ("free-forgetful") functors:
$F: \cat{Sets}\to \cat{Vect}_k $, $ G: \cat{Vect}_k \to \cat{Sets}$.
It is easy to see that for any set $X \in \cat{Sets}$ the unit of adjunction 
$\eta_X: X \to GF (X)$ is given by inclusion of basis, and for any vector space $V$, and counit $\epsilon_V: FG (V) \to V$ is given by continuation of identity linear map on the basis of space $FG(V)$.
\end{ex}

\subsection{Exercises}

\begin{exer}
Suppose we are given three functor $F, G, H: \mathcal C \to \mathcal D$. 
Show that if $\alpha: F \to G$ and $\beta: G \to H$ are natural transformation, then the composition $\beta \alpha: F \to H$ is also natural.
Show that the functors $F: \mathcal C \to \mathcal D$ form a category of $\cat{F}(\mathcal C, \mathcal D)$, where the set of morphisms between any functors is given by natural transformations. Which natural transformation defines identical morphism?
\end{exer}

\begin{exer}
Let $X$ be a topological space. Define a poset $I(X)$ as the set of all closed subsets of $X$ (see Exercise \ref{posetCat} and \ref{posetTop}). and poset $J(X)$ as the set of all subsets of $X$. Show that the functor of embedding of category $ \mathcal C_{I (X) }\hookrightarrow \mathcal C_{J (X) }$ has left adjoint functor that maps an arbitrary subset in $A$ to its closure $\bar{A}$.
\end{exer}

\begin{exer}
Let $I$ be the set of all ideals in commutative ring $\mathbb C [x_1, \dots, x_n]$, and $J$ is the set of all subsets in $\mathbb C^n$ . Define $f: I \to J^{op}$ as a function that takes the ideal to its zero set in $\mathbb C^n$, and $f: J^{op} \to I$
which maps any subset of $\mathbb C^n$ into an ideal of polynomials which annihilates it . 
Show that $f$ and $g$ form Galois correspondence. Conclude that they form a pair of adjoint functors between the respective categories.
\end{exer}

\begin{exer}
Show that the complexifiaction functor $\cat{Vect}_{\mathbb R}\to \cat{Vect}_{\mathbb C}$ (see. Exercise \ref{complexification}) is left to adjoint to functor which restricts the scalars $\cat{Vect}_{\mathbb C}\to \cat{Vect}_{\mathbb R}$.
\end{exer}

\begin{exer}
Let $ G: \cat{Alg}_k \to \cat{Vect}_k $ be forgetful functor. Describe its left adjoint. (\textit{Hint:} Given a vector space $V$ build tensor algebra $T(V) = k \oplus V \oplus (V \otimes V) \oplus \dots $ and show that the correspondence $ V \mapsto T (V) $ determines left adjoint functor to $G$).
\end{exer}

\begin{exer}
Describe the unit and counit of adjunction between $\cat{AbGroup}$ and $ \cat{Group}$.
\end{exer}

\begin{exer}
Describe unit and counit for examples and exercices presented in this section.
\end{exer}

\subsection{References}

\begin{enumerate}
\item Serge Lang, \textit{Algebra}, Springer-Verlag, New York, (2002). 
\item Saunders Mac Lane, \textit{Categories for the working mathematician}, Springer-Verlag, (1971).
\item Masaki Kashiwara, Pierre Schapira, \textit{Categories and Sheaves}, Springer-Verlag, Berlin Heidelberg, 2006.
\item J. Adamek, H. Herrlich, G. Strecker, \textit{Abstract and Concrete Categories: The Joy of Cats}, Wiley-Interscience, (1990).
\item \textit{nLab}(https://ncatlab.org) --- collective and open the online wiki laboratory which contains a lot of useful information on the theory of categories. 

\end{enumerate}

\newpage

\section{Algebra and its radical, connection with quivers}

\subsection{Summary on algebras, their radical, basic algebras.}

Throughout this section we assume that $ A $ is a finite dimensional  algebra over algebraically closed field $ k $, i.e., $A$ is a finite dimensional vector space endowed with an associative multiplication with unit. 

\begin{exs}
\

\begin{itemize}
\item [a)] $ A = k $.
\item [b)] Examples of  infinite dimensional algebra: algebra $ k [x] $ of all polynomials in one variable $ x $ with coefficients in the field $ k $, and algebra  $ k [x_1, \dots, x_n] $ of all polynomials in the variables $ x_1, \dots, x_n $.
\item [c)] Algebra $ k [x] / (x ^ 2) $ of `dual numbers' consists of all pairs of the form $ a + b \ast x $, where $ a, b \in k $, where  $ x $ is an element such that $ x ^ 2 = 0 $. Is easy to see that $ \dim k [x] / (x ^ 2) = 2 $.
\item [d)] if $ A $ is an algebra, then the set $ M_n (A) $ of all $ n \times n $ matrices with coefficients in $ A $ is also an algebra with conventional operations of addition and multiplication of matrices. If  $ A $ is finite dimensional, then $ M_n (A) $ is finite dimensional. In particular, dimension of the algebra $ M_n (k) $ equals $ n ^ 2 $.
\item [e)] Subset 
$$
	\mathbb U_n (k) =\left (\begin{array}{cccc }
	k & k & \dots & k \\
	0 & k & \dots & k \\
	\vdots & \vdots & \dots & \vdots \\
	0 & 0 & \dots & k \\	
\end{array}\right)
$$
of all upper-triangular  matrices $ M_n (k) $ is a subalgebra in $ M_n (k) $.
\item [f)] Associative ring $k \langle x_1, x_2 \rangle$ of all polynomials of two noncommutative variables $ x_1 $ and $x_2$ is a infinite dimensional algebra (called the free algebra on two generators). 
\item [g)] Let $ G $ be a finite group with identity $ e $. The group algebra $ k [G] $ is an  algebra with a basis of $\{a_g\ |\ g \in G \}$ and multiplication $a_ga_h = a_{gh}$. For example, if $ G $ is a cyclic group of order $ m $, then $ k [G] \simeq k [x] / (x ^ m-1) $.
\end{itemize}
\end{exs}

Recall that the \textbf{radical} $ J (A) $  is the intersection of all maximal right ideals in $ A $. It can be shown that the radical $ J (A) $ is the intersection of all maximal left ideals, and therefore $ J (A) $ is two-sided ideal. Algebra $ A $ is called \textbf{semi-simple} if $ J (A) = 0 $. Is easy to see that $ J (A / J (A)) = 0 $.

\begin{exer}\label{alg-basic}
     Let $ A =\mathbb U_2 (K) = \left (\begin{array}{c c}
k & k \\
0 & k
\end{array}\right) $. Show that $$ J (A) =\left (\begin{array}{c c}
0 & k \\
0 & 0
\end{array}\right), $$ and describe radical of algebra $ \mathbb U_n (k) $.
\end{exer}

\begin{exer}\label{mat-semisimple}
	Let $ A = M_{d_1}(k) \oplus \dots \oplus M_{d_n}(k) $ (where $ M_{d_i}(k) $ is a matrix algebra over a field  $k$). Show that $ A $ semi-simple.
\end{exer}

\begin{exer}
Suppose that $\textrm{char} k=0$. Show that the group algebra $ k [G] $ of finite group $ G $ is semi-simple.
\end{exer}

\begin{exer}\label{morRad}
	Let $ f: A \to B $ be a surjective homomorphism of algebras show that $ f (J (A)) = J (B) $.
\end{exer}

Exercise \ref{mat-semisimple} claims that arbitrary algebra  
\begin{equation}\label{al-semisimple}
A = M_{d_1}(k) \oplus \dots \oplus M_{d_n}(k)
\end{equation} 
is semi-simple. In fact, the reverse is also true: theorem of  Joseph Wedderburn  says that an arbitrary semi-simple finite dimensional algebra has the form \eqref{al-semisimple}, and thus for any  algebra $ A $ we have:
\begin{equation}\label{Wudd}
	A / J (A) \simeq M_{d_1}(k) \oplus \dots \oplus M_{d_n}(k).
\end{equation}
We will call the algebra $ A $ \textbf{basic} if in decomposition \eqref{Wudd} all $d_i$  equal $1$, i.e., $ A / J (A) \simeq \prod_i k $. For example, algebra $ A =\mathbb U_n (K) $ is basic since  in this case $ A / J (A) \simeq \underbrace{k \times \dots \times k}_{n\ \text{copies}}$ (see Exercise \eqref{alg-basic}).

\begin{exer}\label{dualN}
Show that the algebra $ A = k [x] / (x ^ m) $  is basic and its radical $ J (A) $  is generated by $ x $.
\end{exer}

\begin{exer}\label{dualMat}
Show that the algebra $ A =\left (\begin{array}{c c}
k & k [x] / x ^ 2 \\
0 & k [x] / x ^ 2
\end{array}\right) $  is basic and its radical has the form $$ J (A) =\left (\begin{array}{c c}
0 & k [x] / x ^ 2 \\
0 & xk [x] / x ^ 2
\end{array}\right). $$
\end{exer}

\subsection{Brief review on representations of algebras}

A \textbf{representation} of an algebra $ A $ (left $ A $-module) is  a vector space $ V $ with algebra homomorphism $ \rho: A \to \textrm{End}(V) $. 

\begin{exs}
\

\begin{itemize}
\item [a)] $ V = 0 $.
\item [b)] $ V = A $, and $ \rho: A \to \textrm{End}A $ is defined as follows: $ \rho (a) $ is an operator of left multiplication by $ a $, i.e., $ \rho (a) b = ab $ (conventional product). Such representation  is called \textbf{regular}.
\end{itemize}
\end{exs}

Given two representations  $ (V_1, \rho_1) $ and $ (V_2, \rho_2) $ of an algebra $ A$, a  \textbf{morphism} between them is defined by the linear operator $ \phi: V_1 \to V_2 $ such that the following diagram commutes
\begin{center}
\begin{tikzcd}[row sep=huge,column sep=huge]
V_1 \arrow[d,"\phi"'] \arrow[r, "\rho_1(a)"] & V_1 \arrow[d,"\phi"] \\
V_2 \arrow[r, "\rho_2(a)"'] & V_2
\end{tikzcd}
\end{center} 
for all $ a \in A $.

Thus, one can create a category  $ \textrm{Rep}A $ of representations of  algebra $ A $. 
Basic algebras play a fundamental role in the theory of representations of finite dimensional  algebras image due to the following theorem

\begin{thm}
	For any finite dimensional algebra $ A $ exists  basic finite dimensional  algebra $ B $ such that categories  $ \textup{Rep}A $ and $ \textup{Rep}B $ are equivalent (see Definition \ref{cat-equiv}).
\end{thm}

So the study of  the representations of all finite dimensional algebra `` reduces '' to the study of the representations  of basic algebras. In the following two subsections we'll make sure that arbitrary basic algebra $ A $ is isomorphic to the quotient algebra of path algebra of  certain quiver.

\subsection{Quivers and path algebras}

A \textbf{quiver} $Q$  is an oriented  graph. We will define quiver $ Q $ by the set of vertices $ Q_0 $, a set  of edges  (arrows) $ Q_1 $, and for a given arrow $ h \in Q_1 $, we will denote by  $ s (h) $, $ t (h) $ its initial and terminal vertex:
\begin{center}
\begin{tikzcd}
s (h) \arrow [r, "h"] & t (h).
\end{tikzcd}
\end{center}

A \textbf{representation} of a quiver  $Q$ is the setting  a vector space $V_i$  for each vertex $ i \in Q_0 $, and   a linear mapping $ V_h: V_{s (h) }\to V_{t (h) }$ for each  arrow $ h \in Q_1 $.

The theory of representations of quivers is closely connected  with the theory of  representations of algebras. Given a quiver $Q$  we will associate an algebra $kQ$ (so-called \textbf{ path algebra} of $Q$ such that representation of the quiver $Q$ are closely connected with the representations of the path algebra $kQ$ (corresponding  categories are equivalent).

The \textbf{path algebra} $kQ$ of a quiver   $ Q $ is an algebra over a field $ k $, the basis of which is formed by  all oriented paths in $ Q $ (including trivial paths  $ p_i, i \in Q_0 $), and multiplication is defined by  concatenation of paths. If  two paths can not be linked, then their product is defined as $0$.

\begin{ex}
The path algebra of a quiver  
\begin{center}
\begin{tikzcd}
1 \arrow [r, "h"] & 2
\end{tikzcd}
\end{center}
has a  basis  of 3 elements  $ p_1 $, $ p_2 $ (trivial paths on the vertices) and $ h $ (path of length 1), endowed  by multiplication $ p_i ^ 2 = p_i, i = 1,2 $, $ p_1p_2 = p_2p_1 = 0 $, $ p_1h = hp_2 = h $, $ hp_1 = p_2h = h ^ 2 = 0 $. So it is  easy to see that there exists  an isomorphism $ kQ \cong \left (\begin{array}{c c}
k & k \\
0 & k
\end{array}\right) $, which is defined by
$$
	p_1 \mapsto \left (\begin{array}{c c}
1 & 0 \\
0 & 0
\end{array}\right), \qquad p_2 \mapsto \left (\begin{array}{c c}
0 & 0 \\
0 & 1
\end{array}\right), \qquad h \mapsto \left (\begin{array}{c c}
0 & 1 \\
0 & 0
\end{array}\right).
$$
\end{ex}
 
\begin{exer}
	Show that the algebra $  kQ $ is generated by  $ p_i, i \in Q_0 $ and $ a_h, h \in Q_1 $ with the following relations:
	\begin{itemize}
		\item [1)] $ p_i ^ 2 = p_i $, $ p_ip_j = 0 $ if $ i \neq j $;
		\item [2)] $ a_hp_{s (h) }= a_h $ $ a_hp_j = 0 $ if $ j \neq s (h) $;
		\item [2)] $ p_{t (h) }a_h = a_h $ $ p_ia_h = 0 $ if $ i \neq t (h) $.
	\end{itemize}
\end{exer}

\begin{exer}\label{path-radical}
Supposing that $Q$ is acyclic, show that the radical of the algebra $kQ$ is generated by all the arrows in $Q$.
\end{exer}

\begin{exer} 
	Using Exercise \ref{path-radical}, show that $ kQ / J (kQ) \simeq \prod_{i \in Q_0}k $. Thus $ kQ $ is a basic algebra.
\end{exer}

\subsection{Quiver of basic algebra.}

Let $A$ be a finite dimensional algebra. Recall that an element of the algebra $ e $ is called \textbf{idempotent} if $ e ^ 2 = e $. Two idempotents $ e, f \in A $ are called \textbf{orthogonal} if $ ef = fe = 0 $. Idempotent $ e $ is called \textbf{primitive} if $ e $ can not be expressed as $ e = e_1 + e_2 $, where $ e_1 $, $ e_2 $ are   nonzero idempotents in $ A $.
Idempotent $ e $ is called  \textbf{central} if $ ae = ea $ for all $ a \in A $.
We say that the algebra $ A $ is \textbf{connected} if $ A $ can not be represented as a direct product of two algebras, or  equivalently  if $ 0 $ and $ 1 $ are  the only central idempotents in $ A $.

Let $ A $ be a basic finite dimensional algebra, and $ \{e_1, \dots, e_n \}$ be a complete set of orthogonal primitive idempotents in $ A $. The \textbf{Gabriel quiver} $Q_A$, associated with $A$, is defined as follows:
\begin{itemize}
	\item ``vertices'' in $ Q_A $ are enumerated by the  elements of set $ \{e_1, \dots, e_n \}$;
	\item Number of ``arrows '' between the vertices  $e_i$ and $e_j$, equals 
	$$ \dim e_i (J (A) / J ^ 2 (A)) e_j. $$
\end{itemize}
One can prove that quiver $Q_A$ does not depend of the choice of the complete set of primitive idempotents in $ A $.

\begin{exs}
\

\begin{itemize}

\item Let $ A = k [x] / (x ^ m) $. It is easy to show that $ e = 1 $ is the  only nonzero idempotent in $ A $. Using Exercise \ref{dualN}, we have  $ J (A) = (x) $, thus $ J ^ 2 (A) = (x ^ 2) $, and therefore $ e \dim (J (A) / J^2(A)) e = 1. $ Therefore, $Q_A$  has the following form:

\
\begin{tikzpicture}
    \node (1) {1} ;
    \Loop[dist=1cm,label=$\alpha$,dir=EA,labelstyle=right](1)  
\end{tikzpicture}

\item Let $ A =\left (\begin{array}{c c}
k & k [x] / x ^ 2 \\
0 & k [x] / x ^ 2
\end{array}\right) $. Is easy to see that
$$ \left \{e_1 = \left (\begin{array}{c c}
1 & 0 \\
0 & 0
\end{array}\right), e_2 = \left (\begin{array}{c c}
0 & 0 \\
0 & 1
\end{array}\right) \right \}$$
is a complete set of primitive orthogonal idempotents (therefore $Q_A$ has two vertex), and 
$$ J (A) =\left (\begin{array}{c c}
0 & k [x] / x ^ 2 \\
0 & xk [x] / x ^ 2
\end{array}\right), \quad J ^ 2 (A) =\left (\begin{array}{c c}
0 & xk [x] / x ^ 2 \\
0 & 0
\end{array}\right) $$
(See Exercise \ref{dualMat}). So $ \dim J (A) / J ^ 2 (A) = 2 $, therefore  $Q_A$  has two arrows.
Direct calculations give the following 
$$
	\dim e_i (J (A) / J ^ 2 (A)) e_j =\left \{\begin{array}{l}1, \quad (i, j) = (1, 2) \mbox{or }(i, j) = (2,2), \\0, \quad \mbox{otherwise.}
	\end{array}\right.
$$

\end{itemize}
\end{exs}

If $ Q $ is an arbitrary  finite quiver,  let $ R_Q $ be two-sided ideal in algebra $kQ$  generated  by arrows of $ Q $. We say that two-sided ideal $ I \subset kQ $ \textbf{admissible} if
$$
	R_Q ^ m \subseteq I \subseteq R_Q ^ 2,
$$
for a $ m \geq 2 $. In other words, $ I $ is admissible, if it does not  contain  arrows  of $ Q $ and includes all paths of length $ \geq m $. If $ I $ is admissible, the quotient  algebra $ kQ / I $ is called  \textbf{bound path  algebra}.

The following theorem defines the canonical form of basic finite dimensional algebras.

\begin{thm}
Let $ A $ be basic, connected, finite dimensional  algebra. There exists  an admissible  ideal $ I $ in $ kQ_A $ such that $ A \cong kQ_A / I $.
\end{thm}

\begin{exer}
	Build  the Gabriel quiver of the algebra $ A =\mathbb U_n (k) $.
\end{exer}

\begin{exer}
	Let $ A =\mathbb U_3 (k) $, and $ C $ be  subalgebra consisting of all matrices 
	$$
		\lambda =\left (\begin{array}{cc c}\lambda_{11}& \lambda_{12}& \lambda_{13}\\
		0 & \lambda_{22}& \lambda_{23}\\
	    0 & 0 & \lambda_{33}\end{array}\right)
	$$
	such that  $ \lambda_{11}=\lambda_{22}=\lambda_{33}$. Show that $ C $ is isomorphic to $ kQ / I $, where $ I =\langle \alpha ^ 2, \beta ^ 2, \alpha \beta \rangle $  is ideal in $ kQ $, and $ Q $ is the  following quiver 
\
\begin{center}
\begin{tikzpicture}
    \node (1) {1};
	\Loop[dist=1.5cm,label=$\alpha$,dir=WE,labelstyle=left](1)  
	\Loop[dist=1.5cm,label=$\beta$,dir=EA,labelstyle=right](1)  
\end{tikzpicture}
\end{center}
\end{exer}
\subsection{References}

\begin{enumerate}
\item Y.A. Drozd, V.V. Kirichenko, \textit{Finite dimensional algebras}, Springer-Verlag, Berlin/New York (1995).
\item Gabriel P., Roiter A.D., \textit{Representations of Finite-Dimensional Algebras}, Springer (1997).
\item Assem I., Simson D., Skowronski A., \textit{Elements of the representation theory of associative algebras}, vol 1, London Mathematical Society Student Texts, (2007). 
\item Auslander M., Reiten I., Smalo SO, \textit{Representation theory of Artin algebras}, Cambridge University Press, (1977).
\item Etingof P., Golberg O., Hensel S., Liu T. \textit{Introduction to Representation Theory}, American Mathematical Society, (2011).
\item \textit{fdLit}(http://www.math.uni-bonn.de/people/schroer/fd-literature.html) --- Jan Schr\"{o}er's page containing literature on finite dimensional algebras, their representations and related topics. 
\end{enumerate}
 
\newpage

\section{Path algebra as a left adjoint functor}

\epigraph{``Adjoint functors arise everywhere.'' }{--- \textup{Saunders Mac Lane}}

In the previous lesson we saw that quivers (and their representations) play a fundamental role in the structure and representation theory of finite dimensional algebras. In this lecture we will try to show that the construction ``quiver'' $ \leftrightarrows $ ``algebra'' ($Q \mapsto kQ$, $A \mapsto Q_A$) can be interpreted as a pair of adjoint functors between certain categories. This lecture is based on a part of a joint research project with John William Macquarrie (UFMG).

Consider the category $\cat{Quiv}$, where the objects are finite quivers and the morphisms are embeddings of quivers (i.e., embeddings on vertices and on corresponding sets of arrows). 
Let $\cat{Quiv}^{ac}$ be the category of acyclic quivers.
Consider the category $ \cat{SBAlg}$, where  objects are finite dimensional basic algebras and morphisms are surjective algebra homomorphisms. 

\begin{exer}\label{cCatkFunct}
Show that $\cat{Quiv}$ and $ \cat{SBAlg}$ are categories, $\cat{Quiv}^{ac}$ is a full subcategory in $\cat{Quiv}$ (i.e., morphisms between two objects in $\cat{Quiv}^{ac}$ are the same as in $\cat{Quiv}$). Show that the correspondence $Q \mapsto kQ$ generates the contravariant functor 
$K[-]: \cat{Quiv}^{ac}\to \cat{SBAlg}$.
\end{exer}

Is easy to see that the correspondence $A \mapsto Q_A$ does not generate a functor between categories $\cat{SBAlg}$ and $ \cat{Quiv}$. Indeed, the choice of a complete set of primitive orthogonal idempotents in $A$, in general, is not unique, and the choice of the basis in the space $e(J (A)/J^2(A))f$ (arrows between $e$ and $f$) is not canonical. 
We consider certain intermediate category between ``quivers'' and ``algebras'' so that the correspondence above defines a functor.

\subsection{Quotient category $\cat{SBAlg}_n$}

The construction of quotient category is similar to the construction of quotient set or quotient algebra. Let  $\mathcal C$ be an arbitrary category. Assume that the equivalence relation $\sim$ is defined on morphisms $\cat{Mor}(\mathcal C)$. That is, for arbitrary $X, Y \in \mathcal C$ the set $ \cat{Mor}(X, Y) $ splits  into equivalence classes $[\alpha]$, which satisfy the condition: once $[\alpha] = [\alpha ']$, then $[\beta \alpha] = [\beta \alpha']$ and $ [\alpha \beta] = [\alpha '\beta] $, when the composition of morphisms makes sense. Now we can create a new category of $\mathcal C / \sim $, which will be called \textbf{quotient category}. Objects in $\mathcal C /\sim$ are the same as the objects in 
$\mathcal C$, and morphisms set $\cat{Mor}_{\mathcal C / \sim}(X, Y)$ consists of equivalence classes of morphisms $ \cat{Mor}_{\mathcal C}(X, Y)$ relative to $\sim$. Composition of morphisms is given by the rule $[\beta][\alpha]=[\beta \alpha]$.

Define the following equivalence relation in the category $\cat{SBAlg}$. Let $A, B \in \cat{SBAlg}$ and $\alpha_1, \alpha_2 \in \cat{Mor}(A, B) $. We say that $\alpha_1$ and 
$\alpha_2$ are \textit{$n$-depth} (denoting this by $\alpha_1 \sim_n \alpha_2 $) if 
$$
	(\alpha_1 - \alpha_2) (J ^ i (A)) \subseteq J ^{i + 1}(B), \qquad 0 \leq i \leq n,
$$
setting $ J ^ 0 (A) = A $.

\begin{exer}
Show that $\sim_n$ is an equivalence relation on $\cat{SBAlg}$. (\textit{Hint:} use Exercise \ref{morRad}).
\end{exer}

Thus, we form quotient category $$ \cat{SBAlg}_n =\cat{SBAlg}/ \sim_n.$$ Denote by $\Pi_n: \cat{SBAlg}\to \cat{SBAlg}/ \sim_n $ the corresponding quotient functor (which is identical on objects, and maps each morphism $\alpha: A \to B$ to its equivalence class $ [\alpha] _n $ towards $ \sim_n $).

\subsection{Category of Vquivers}

Finite (pointed) \textbf{Vquiver} $$VQ = (VQ_0, VQ_1)$$ is given by a finite set of vertices 
$VQ_0^* = \{* \}\cup VQ_0 $, where $ VQ_0 = \{e_1, \hdots, e_n \}$, with finite-dimensional vector space $VQ_{e, f}$ for every pair of vertices $e, f \in VQ_0^*$ such that $VQ_{*,e}= VQ_{e, * }= 0$ for all $ e $. 

Denote by $\Sigma_{VQ}$ free $k$-module generated by $VQ_0$, which we treat as a semi-simple algebra, defining 
$$ e_i \cdot e_j = \left \{\begin{array}{l l}1, \quad & i = j, \\0 & i \neq j  \end{array}\right..$$ 
By $VQ_1$ we denote the direct sum  $\bigoplus_{e, f \in VQ_0}VQ_{e, f}$, which has the natural structure of $\Sigma_{VQ}$-bimodule.

\begin{defn}
A \emph{map of finite Vquivers} $\rho:VQ \to VR$ consists of 
\begin{itemize}
\item a pointed map $\rho_0:VQ_0^*\to VR_0^*$ (i.e., such that $\rho_0(*)=*$) that restricts to a bijection from the elements of $VQ_0$ not mapping to $*$, onto $VR_0$.
\item a linear map $\rho_{e,f}:VQ_{e,f}\to VR_{\rho_0(e),\rho_0(f)}$ for each pair of vertices $e,f\in VQ_0^*$.
\end{itemize} 

We say that $\rho$ is \emph{surjective} if every $\rho_{e,f}$ is surjective.
\end{defn} 

\begin{rem}
One can check that $\rho$ is a morphism $\rho:VQ \to VR$ iff $\rho_0$ is an algebra homomorphisms and $\rho_{VQ}=\bigoplus_{e,f\in VQ_0} \rho_{e,f}:VQ_1\rightarrow VR_1$ is a $\Sigma_{VQ}-\Sigma_{VQ}$ bimodule morphism. 
\end{rem} 

Now define the category $\cat{SVQuiv}$ which has objects finite dimensional Vquivers and morphisms surjective maps of Vquivers.  
We say that Vquiver $ VQ = (VQ_0, VQ_1) $ is \textbf{acyclic} if there is $n>0 $ such that
$$
	\underbrace{VQ_1 \otimes_{\Sigma}\dots \otimes_{\Sigma}VQ_1}_{n}= 0.
$$
Now denote by $\cat{SVQuiv}^{ac}$ the full subcategory consisting of acyclic Vquivers.

\begin{exer}
Define the natural contravariant functor between categories $\cat{Quiv}$ and $\cat{SVQuiv}$.
\end{exer}

\subsection{Functor ``path algebra''.} Given a pair $(\Sigma, V)$ (where $\Sigma$ is an arbitrary algebra and $V$ is a $\Sigma-\Sigma$-bimodule) associate the \textbf{tensor algebra} $T(\Sigma, V)$ by
$$
	T (\Sigma, V) =\Sigma \oplus V \oplus V \otimes _ \Sigma V \oplus \dots.
$$
The tensor algebra satisfies the following universal property.
\begin{prop}\label{UP of tensor algebra}
Let $A$ be a $k$-algebra.  Let $\Sigma$ be a semi-simple $k$-algebra and $V$ be a $\Sigma-\Sigma$ bimodule. Suppose we have two functions 
$$\varphi_0:\Sigma\to A, \,\,\varphi_1:V\to A$$ 
such that
\begin{enumerate}
\item $\varphi_0$ is a $k$-algebra homomorphism
\item $\varphi_1$ is a $\Sigma-\Sigma$ bimodule homomorphism, with $A$ treated as a $\Sigma-\Sigma$ bimodule via $\varphi_0$.
\end{enumerate}
Then there is a unique $k$-algebra homomorphism $\varphi:T(\Sigma,V)\to A$ such that $\varphi|_{\Sigma} = \varphi_0, \varphi|_V = \varphi_1$. 
\end{prop}

Let $VQ = (VQ_0 ^ *, VQ_{e, f})$ be a finite acyclic Vquiver. 
The \emph{path algebra} $k[VQ]$ is the tensor algebra $T(\Sigma_{VQ},VQ_1)$ defined above. We have that $k[VQ]$ is a basic algebra, therefore $k[VQ]\in \cat{SBAlg}$. 
Let $\rho : VQ \to VR$ be a surjective map of Vquivers.  
We obtain the maps $\varphi_0, \varphi_1$ as in the previous proposition in the obvious way. Indeed we have that $k[VR]=T(\Sigma_{VR},VR_1)$. 
As $\rho_0$ is surjective on vertices therefore it generates a surjective algebra morphism 
$\varphi_0:\Sigma_{VQ} \twoheadrightarrow \Sigma_{VR}\subset k[VR]$. 
On the other hand $\rho_1$ generates a map $$\varphi_1:VQ_1\rightarrow VR_1\subset k[VR].$$ 

\begin{exer}
Check that $\varphi_1$ is a $\Sigma_{VQ}-\Sigma_{VQ}$-bimodule homomorphism. 
\end{exer} 

Hence we obtain a unique $k$-algebra homomorphism $k[\rho]:k[VQ]\to k[VR]$ by Proposition \ref{UP of tensor algebra}. As $\rho$ is surjective map one shows (check!) that $k[\rho]$ is a surjective algebra homomorphism.

Thus, we get:
\begin{prop}
The following construction defines a covariant functor:
\begin{align*}
k [-]: & \cat{SVquiv}^{ac}\rightarrow \cat{SBAlg}, \\
\mathscr{K}_n [-] =\Pi_n \circ k [-]: & \cat{SVquiv}^{ac}\rightarrow \cat{SBAlg}_n.
\end{align*}
\end{prop}

\begin{rem}
Note that the composition of the functors $ k [-] \circ V [-] $ is a functor from Exercise \ref{cCatkFunct}.
\end{rem}

\subsection{Functor ``Gabriel Quiver''}

Now we define the functor in another direction. Let $A$ be a finite dimensional algebra. Recall Wedderburn-Malcev theorem (in a short form):

\begin{thm}
There is a subalgebra $\Sigma$ of $A$ such that $A=\Sigma\oplus J(A)$ as $k$-vector spaces and $\Sigma\cong A/J(A)$ as algebras. 
For two any subalgebras $\Sigma$ e $\Sigma'$ such that $A=\Sigma\oplus J(A)=\Sigma'\oplus J(A)$, there exists $w\in J(A)$ such that 
    $$
    	\Sigma'=(1+w)\Sigma(1+w)^{-1}.
    $$
\end{thm}

The key idea is to define the ``vertices'' in Gabriel quiver as the orbit under $J(A)$.
Let $A \in \cat{SBAlg}$. Denote by $ \pi_A: A \twoheadrightarrow A / J (A) $ the canonical projection. Algebra homomorphism $s: A /J(A) \hookrightarrow A$ is called  \textit{spllitting}of $\pi_A$, if $ \pi_A \circ s =\textrm{id}_{J(A)}$. By $ \mathcal S_A $ denote the set of all such splittings and by $ \mathcal E_A $ the collection of all possible complete sets  primitive orthogonal idempotents in $A$. Recall that, according to the theorem Wedderburn-Malcev, $\mathcal S_A$ is non-empty. Since $A$ is basic then $A/J(A) \cong \prod_{i = 1}^{n}k $. So any split  
$s \in S_A$ defines a complete set of orthogonal primitive idempotents
$ \{s (j_1), \dots, s (j_n) \}\in \mathcal E_A $, in which $j_1, \dots, j_n$ is the unique  complete set of primitive orthogonal idempotents $A/J(A)$. 
Denote the above correspondence by $\Phi:\mathcal S_A\rightarrow \mathcal E_A$. 

\begin{exer}
Show that $ \Phi $ is a bijection.
\end{exer}

Any element $w\in J(A)$ gives rise to an automorphism 
$$
	a\mapsto (1+w)a(1+w)^{-1},\qquad a\in A.
$$
We denote such automorphism by ${}^{1+w}(-)$, respectively we use the notation ${}^{1+w}a:=(1+w)a(1+w)^{-1}$. Let $\mathcal G(A) \triangleleft \textrm{InnAut}(A)$ be the group of all such automorphisms. We will use the notation $\mathcal G$ in cases when it is clear which algebra we are using. Denote by ${}^{\mathcal G}a = \{{}^{1+w}a\,|\,w\in J\}$ the orbit 
of a given element $a\in A$.

Let $A$ be a basic finite dimensional algebra and $s\in \mathcal S_A$ be any splitting, and let $\Phi(s)\in \mathcal E_A$ be the corresponding set of primitive orthogonal idempotents in $A$. Define the \textit{Vquiver} $GQ(A)$ of $A$ as follows:
\begin{align*}
	GQ(A)_0&:=\{*\}\cup \{{}^{\mathcal G}e\,|\,e\in \Phi(s)\},\\
    GQ(A)_{{}^{\mathcal G}e,{}^{\mathcal G}f}& := e\frac{J(A)}{J^2(A)}f,\qquad \mbox{for fixed}\ \ e,f\in \Phi(s).
\end{align*}

\begin{exer}
Show that \textit{Vquiver} $GQ(A)$ is well-defined.  That is, it does not depend on the choice of $s\in \mathcal S_A$.  
\end{exer}

Let $A,B$ be algebras with inner automorphism groups $\mathcal G = \mathcal G(A)$ and $\mathcal H = \mathcal G(B)$.
Given  a surjective algebra homomorphism $\alpha:A \to B$, we define the map of Vquivers $GQ(\alpha):GQ(A)\to GQ(B)$ as follows:
\begin{align*}
GQ(\alpha)({}^{\mathcal G}e)= &{}^{\mathcal H}\alpha(e);\\
GQ(\alpha):&e\dfrac{J(A)}{J^2(A)}f\to \alpha(e)\dfrac{J(B)}{J^2(B)}\alpha(f)\\
&e(j+J^2(A))f\mapsto \alpha(e)(\alpha(j)+J^2(B))\alpha(f). 
\end{align*} 

\begin{exer}
Show that the map in the definition above is well-defined map of Vquivers.
Moreover $GQ(\alpha)$ is surjective. 
\end{exer}

Thus, we have
\begin{prop}
The construction above defines a covariant functor
$$
	GQ (-): \cat{SBAlg}\rightarrow \cat{SVquiv}.
$$
\end{prop}

\begin{exer}\label{continuF}
For any $n\geq 1$ show that there is a unique functor  
$$
	\mathscr{GQ}_n (-): \cat{SBAlg}_n \rightarrow \cat{SVquiv}
$$
such that the following diagram commutes
\begin{center}
\begin{tikzcd}[row sep = huge, column sep = huge]
\cat{SBAlg}\arrow [d, "\Pi_n" '] \arrow [r, "GQ (-)"] & \cat{SVquiv}\\
\cat{SBAlg}_n \arrow [ru, "\mathscr{GQ}_n (-)" '] & 
\end{tikzcd}
\end{center}
Build functor $ \mathscr{GQ}_n (-) $.
\end{exer}

\subsection{Adjunction between functors.} 

By $\cat{SBAlg}^{ac}$ denote the full subcategory of basic algebra such that Vquiver $GQ(A)$ is acyclic. 
Summing up the construction described above, we have the following diagram

\

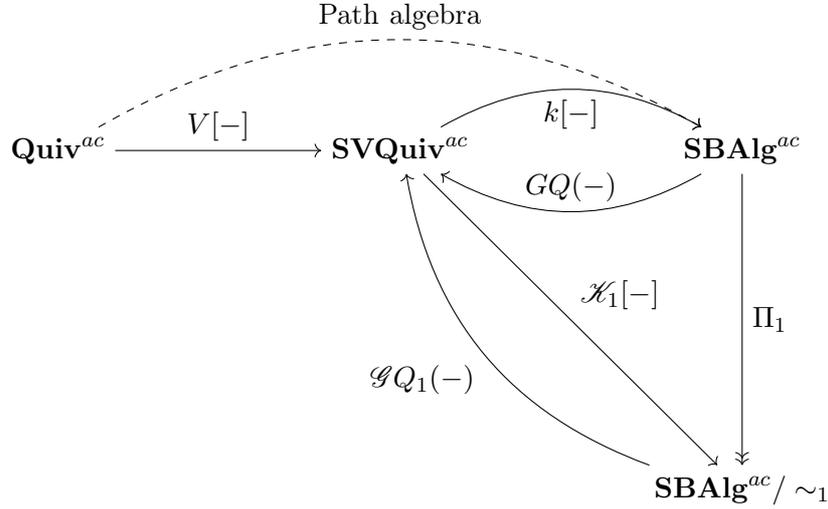
\begin{center}
\begin{tikzpicture}[node distance = 4.5cm, auto]
  \node (A){$ \cat{Quiv}^{ac}$ };
  \node (B) [right of = A]{$ \cat{SVQuiv}^{ac}$ };
  \node (C) [right of = B]{$ \cat{SBAlg}^{ac}$ };
  \node (D) [below of = C]{$ \cat{SBAlg}^{ac}/ \sim_1 $ };
  \draw [->] (A) to node{$ V [-] $ }(B);
  \draw [->, bend left] (B) to node [swap]{$ k [-] $ }(C);
  \draw [->, bend left] (C) to node [swap]{$ GQ (-) $ }(B);
  \draw [->, dashed, bend left] (A) to node{Path algebra}(C);
  \draw [- >>] (C) to node{$ \Pi_1 $ }(D);
  \draw [->] (B) to node{$ \mathscr K_1 [-] $ }(D);
  \draw [->, bend left] (D) to node{$ \mathscr GQ_1 (-) $ }(B);
\end{tikzpicture}
\end{center}
Here the functor $\mathscr K_1 [-]$ denotes the composition of functors $\Pi_1 \circ k [-]$. The functor $\mathscr GQ_1(-)$ is the restriction of the functor constructed in Exercise \ref{continuF} to the subcategory $\cat{SBAlg}^{ac}$. 

The following theorem holds.

\begin{thm}\label{MainAdjunction}
	The functor $\mathscr K_1 [-]$ is left adjoint to $\mathscr GQ_1 (-)$.
\end{thm}

\begin{exer}
As a consequence of the previous theorem, show that an arbitrary algebra $A \in \cat{SBAlg}$ is a quotient algebra of path algebra.
\end{exer}

\begin{exer}
Describe unit and counit of adjunction in Theorem \ref{MainAdjunction}. Describe the image functor $ \mathscr K_1 [-] $. This functor full and faithfull?
\end{exer}

\section*{Acknowledgements}
The first version of these lectures were given on XI summer school ``Algebra, Topology, Analysis'' (1--14 August 2016, Odessa, Ukraine). Authors thanks Sergey Maksimenko for invitation and the presentations of this course. Author also grateful to Nataliia Goloshchapova and Volodymyr Tesko for their helpful remarks, who attentively proofread the text and gave various helpful comments.
The actual version of the course was given on Summer School of Algebra at the Federal University of Minas Gerais (01--03 February, 2017, Belo Horizonte, Brasil). Author thanks John Macquarrie for the kind invitation and for numbers of helpful remarks about the text.

\end{document}